\documentclass[12pt,russian]{article}

\usepackage[T2A]{fontenc}
\usepackage{amsmath,amsfonts,amssymb,amsthm}
\usepackage[russian,english]{babel}

\renewcommand{\No}{N}

\newtheorem{Theorem}{Теорема}
\newtheorem{Lemma}{Лемма}

\theoremstyle{remark}
\newtheorem{Note}{Замечание}
\theoremstyle{plain}

\newenvironment{Proof}{\par\textbf{Доказательство.} }{\hfill\rule{0.5em}{0.5em}}
\newenvironment{VarProof}{\par\textbf{Доказательство}}{\hfill\rule{0.5em}{0.5em}}

\newcommand{\PutItem}[3]{{\textrm{#1} \textit{#2}// \textrm{#3}}}

\newcommand{\PutItemAa}[6]{\PutItem{#1}{#2}{#3 #4. \textbf{#5}. #6.}}

\newcommand{\PutItemAb}[7]{\PutItem{#1}{#2}{#3 #4. \textbf{#5}, \No{}~#6. #7.}}

\newcommand{\PutItemAd}[5]{\PutItem{#1}{#2}{#3 #4. #5.}}

\newcommand{\PutItemBa}[5]{\PutItem{#1}{#2}{#3: #4 #5.}}

\newcommand{\PutItemBb}[4]{\PutItem{#1}{#2}{#3. #4.}}

\newcommand{\PutItemNAb}[6]{\PutItem{#1}{#2}{#3 #4. \textbf{#5}, \No{}~#6.}}

\newcommand{\GE}{\geqslant}
\newcommand{\LE}{\leqslant}

\newcommand{\Ai}{\mathop\mathrm{Ai}\nolimits}
\newcommand{\Bi}{\mathop\mathrm{Bi}\nolimits}

\newcommand{\RE}{\mathop\mathrm{Re}\nolimits}
\newcommand{\IM}{\mathop\mathrm{Im}\nolimits}

\newcommand{\OB}[1]{O\left(#1\right)}
\newcommand{\eps}{\varepsilon}

\newcommand{\node}{-\frac{i}{\sqrt{3}}}

\newcommand{\lk}{<<}
\newcommand{\pk}{>>}

\begin{document}
%%%%%%%%%%%%%%%%%%%%%%%%%% English Title %%%%%%%%%%%%%%%%%%%%%%%%%%%%%%%%%

\begin{center}

{\LARGE\bf
On the Orr--Sommerfeld Equation \\[0.3cm]
with Linear Profile
}

\vspace{10mm}

{\large
A. V. Dyachenko$^{(a)}$ and A. A. Shkalikov$^{(b)}$
}

\vspace{10mm}

Department of Mechanics and Mathematics,\\
Moscow State University,\\
Moscow, Russia\\
{\em
${}^{(a)}$e-mail: sasha\_d@auriga.ru\\
${}^{(b)}$e-mail: ashkalikov@yahoo.com\\
}

\vspace{15mm}

{\bf Abstract}

\end{center}

\noindent
The Orr--Sommerfeld equation with linear profile on the finite interval
is considered. The behavior of the spectrum of this problem is completely
investigated for large Reynolds numbers. The limit curves are found to which
the eigenvalues concentrate and the counting eigenvalue functions along
these curves are obtained.
\vskip 2mm

\newpage

%%%%%%%%%%%%%%%%%%%%%%%%%% Russian Title %%%%%%%%%%%%%%%%%%%%%%%%%%%%%%%%%

\begin{center}

{\LARGE\bf
Уравнение Орра--Зоммерфельда \\[0.3cm]
с линейным профилем
}

\vspace{10mm}

{\large
А. В. Дьяченко$^{(a)}$ и А. А. Шкаликов$^{(b)}$
}

\vspace{10mm}

Механико-математический факультет,\\
Московский Государственный Университет,\\
Москва, Россия\\
{\em
${}^{(a)}$e-mail: sasha\_d@auriga.ru\\
${}^{(b)}$e-mail: ashkalikov@yahoo.com\\
}

\vspace{15mm}

{\bf Аннотация}

\end{center}

\noindent
Рассматривается уравнение Орра--Зоммерфельда с линейным профилем. Полностью
изучено поведение спектра для больших чисел Рейнольдса. Найдены предельные
кривые, около которых концентрируются собственные значения, и получены
функции распределения собственных значений вдоль этих кривых.

\newpage

\section*{Введение}
В статье рассматривается спектральная задача Орра--Зоммерфельда
\begin{gather}
\{(D^2-\alpha^2)^2-i\alpha
R[q(x)(D^2-\alpha^2)-q''(x)]\}y=-\lambda(D^2-\alpha^2)y, \label{orsom_gen}\\
y(-1)=y'(-1)=y(1)=y'(1)=0. \label{orsom_gen_bound}
\end{gather}
Здесь $D=d/dx$, $\alpha$ --- волновое число, $R$ --- число Рейнольдса,
характеризующее вязкость жидкости, а $q(x)$ --- профиль скорости течения
жидкости в канале $|x|<1$. Эта задача получается после линеаризации уравнений
Навье--Стокса для плоскопараллельных течений между двумя фиксированными
стенками.

Свойства задачи Орра--Зоммерфельда во многом определяются свойствами
модельной задачи вида
\begin{gather}
-\eps y''+q(x)y=\lambda y, \label{model_gen}\\
y(-1)=y(1)=0, \label{model_gen_bound}
\end{gather}
где $\eps=1/i\alpha R$ --- малый параметр, лежащий на отрицательной мнимой полуоси.

Хорошо известно, что спектр задачи Орра--Зоммерфельда на конечном интервале
дискретен. Важной является задача изучения поведения собственных значений при
больших числах Рейнольдса $R$, что соответствует малой вязкости жидкости.
Наибольший интерес представляют два стационарных профиля скорости: $q(x)=x$ и
$q(x)=x^2$. Первый называется профилем Куэтта, второй --- профилем Пуазейля.

Задача Орра--Зоммерфельда изучалась многими авторами. Основные результаты и
литературные ссылки можно найти в обзоре Редди, Хеннингсона и Шмидта
\cite{Reddy_Schmidt_Henningson_1993}, монографиях Драйзина и Рида
\cite{Drazin_Reid_1981}, Дикого \cite{Dikiy_1973}, а также в работах
Гейзенберга, Вазова, Лина и др. (см. библиографию в \cite{Drazin_Reid_1981}).

Однако, описание портрета поведения собственных значений этой задачи  при
$\eps\to 0$ (т. е. $R\to\infty$) не было приведено полностью. В этой связи
укажем важные работы Моравец \cite{Morawetz_1952} и Чапмана
\cite{Chapman_2002}. Моравец показала, что при $q(x)=x$ собственные значения
задачи Орра--Зоммерфельда могут локализоваться только вблизи отрезков
$[-1,-i/\sqrt{3}]$, $[1,-i/\sqrt{3}]$ и луча $[-i/\sqrt{3},-i\infty)$, хотя в
\cite{Morawetz_1952} подчеркивается, что не удается получить информацию о
собственных значениях в малых окрестностях первых двух отрезков. Для
$q(x)=x^2$ предположения работы \cite{Morawetz_1952} не реализуются.

Компьютерные программы для вычисления собственных значений задачи
Орра--Зоммерфельда реализовывались Редди, Хеннингсоном и Шмидтом
\cite{Reddy_Schmidt_Henningson_1993}, Трефезеном \cite{Trefethen_1996},
Шкаликовым и Нейманом-заде \cite{Neyman_Shkalikov_2002}, Чапманом
\cite{Chapman_2002}. В двух последних работах было понятно, что в отличие от
модельной задачи, собственные значения задачи Орра--Зоммерфельда локализуются
вдоль двух линий снизу и сверху от отрезков $[\pm 1,\node]$, причем эти линии
при $\eps\to 0$ сливаются в эти отрезки. В работе \cite{Chapman_2002}
приводились некоторые объяснения этого явления, но аналитической формы линий
и явных формул для собственных значений на этих линиях (или функций
распределения) выписано не было. Все это: форма линий и формулы для
собственных значений, --- будет найдено в этой статье.

Необходимо отметить, что многие авторы рассматривали задачу
Орра--Зоммерфельда в связи с вопросом об устойчивости течения жидкости, что
эквивалентно отсутствию собственных значений в верхней полуплоскости.
Устойчивость плоскопараллельного течения Куэтта \lk почти\pk{} доказана в
работе \cite{Romanov_1973}. При $R\to\infty$ в ней используется аппарат
специальных функций (функций Эйри), а при малых значениях числа Рейнольдса
применяются численные расчеты на компьютере. Вопрос об устойчивости для
уравнения Орра--Зоммерфельда рассмотрен и в книге \cite{Zhuk_2001}, где
используются асимптотические методы для нахождения решений.

Модельная задача \eqref{model_gen}, \eqref{model_gen_bound} с линейным
профилем рассматривалась в работах Трефезена \cite{Trefethen_1996} и Редди,
Хеннингсона, Шмидта \cite{Reddy_Schmidt_Henningson_1993}. Аналитическое
объяснение портрета собственных значений этой задачи при $\eps\to 0$ было
проведено в \cite{Shkalikov_1997}, а именно, было доказано, что при $q(x)=x$
собственные значения модельной задачи \eqref{model_gen},
\eqref{model_gen_bound} локализуются вдоль луча $[-i/\sqrt{3},-i\infty)$ и
двух отрезков $[-1,-i/\sqrt{3}]$, $[1,-i/\sqrt{3}]$, а также найдена
асимптотика собственных значений в окрестности указанных отрезков. Более
частный результат независимо получен в \cite{Stepin_1998}.

Итак, рассматривается спектральная задача Орра--Зоммерфельда с линейным 
профилем (течение Куэтта)
\begin{align}
-&i\eps(w''-\alpha^2 w)=(x+i\eps\tilde\lambda)w,\qquad w=y''-\alpha^2 y,
\label{orsom1}\\
&y(\pm 1)=y'(\pm 1)=0. \label{orsom2}
\end{align}
Здесь $\eps=1/\alpha R$, $\tilde\lambda$ --- спектральный параметр, $R$ ---
число Рейнольдса, а $\alpha$ --- волновое число. Для наших рассмотрений
удобно ввести другой спектральный параметр
$\lambda=i\eps(\alpha^2-\tilde\lambda)$.

Соответствующая модельная задача имеет вид:
\begin{align}
-&i\eps y''=(x-\lambda)y,
\label{eq_model1}\\
&y(-1)=y(1)=0, \label{eq_model2}
\end{align}
где $\lambda$ --- спектральный параметр, а $\eps>0$ --- малый параметр.

\begin{figure}[tbp]
\unitlength=60mm
\begin{picture}(2.2,2.1)(-1.4,-2)
\input{pic_model/base.tex}
\end{picture}
\caption{\small Спектр модельной задачи для профиля Куэтта, $R=3000$,
$\alpha=1$} \label{pic_tie}
\end{figure}

Наша цель~--- описать асимптотическое поведение собственных значений
модельной задачи при $\eps\to+0$ на мнимой оси, уточнить асимптотики
\cite{Shkalikov_1997} на отрезках $[\pm 1,\node]$ и оценить число собственных
значений вблизи узла. Затем, пользуясь наработанными методами, сделать то же
самое для задачи Орра--Зоммерфельда с линейным профилем.

План статьи следующий.

В параграфе \ref{par31} рассмотрена модельная задача \eqref{eq_model1},
\eqref{eq_model2}. В теореме \ref{th_orsommod_ray} вычислена асимптотика
собственных значений на луче $[\node,-i\infty)$. В теореме
\ref{th_orsommod_seg} сделано уточнение асимптотики \cite{Shkalikov_1997} на
отрезках $[\pm 1,\node]$ с учетом леммы \ref{lm_aias}. В теореме
\ref{th_orsommod_Nlam} дается удобное представление для функции распределения
собственных значений вне узловой точки $\node$ (в которой не работают
асимптотические приближения). Наконец, оценка количества собственных значений
в круге малого радиуса с центром в узловой точке $\node$ дана в теореме
\ref{th_orsommod_Ndel}. Отметим, что уменьшение радиуса этого круга вплоть до
$\eps^{1/2}\ln\eps^{-\sigma}$ стало возможным благодаря лемме \ref{lm_aias},
в которой расширена стандартная область действия классических асимптотик
функций Эйри.

В параграфе \ref{par32} найдены четыре независимых решения уравнения
Орра--Зоммерфельда с линейным профилем и выписано характеристическое
уравнение краевой задачи в виде определителя размером $2\times 2$, состоящего
из интегралов от функций Эйри.

В параграфе \ref{par33} изучается поведение собственных значений $\lambda$ в
окрестности отрезка $[-1,\node]$. Лемма \ref{lm_uxi_prop} дает необходимые
асимптотические представления для интегралов, входящих в характеристический
определитель. Дальнейший анализ показывает, что собственные значения
локализуются не вдоль самого отрезка $[-1,\node]$, а вдоль двух кривых,
находящихся по разные стороны от этого отрезка и отстоящих от него на
расстояние порядка $\eps^{1/2}|\ln\eps|$. Этот результат сформулирован в
теореме \ref{th_orsom_main}, где также даны асимптотические представления для
этих кривых и самих собственных значений.

Наконец, в параграфе \ref{par34} мы покажем, что в окрестности луча
$[\node,-i\infty)$ собственные значения $\lambda$ ведут себя точно так же,
как и в модельной задаче. Соответствующие формулы выписаны в теореме
\ref{th_orsom_ray}.

\begin{figure}[tbp]
\unitlength=60mm
\begin{picture}(2.2,2.1)(-1.4,-2)
\input{pic_orsom/base.tex}
\end{picture}
\caption{\small Спектр задачи Орра--Зоммерфельда для профиля Куэтта, $R=3000$,
$\alpha=1$} \label{pic_tie_orsom}
\end{figure}

\section{Поведение спектра модельной задачи}
\label{par31}

Замена $\xi=(-i\eps)^{-1/3}(x-\lambda)$ приводит задачу \eqref{eq_model1},
\eqref{eq_model2} к граничной задаче для уравнения Эйри
\begin{align*}
&z''(\xi)=\xi z(\xi),\\
&z(\xi_1)=z(\xi_2)=0,
\end{align*}
где $\;\xi_1=(-i\eps)^{-1/3}(-1-\lambda),\quad\xi_2=(-i\eps)^{-1/3}(1-\lambda)$.

Известно \cite{Olver_1974}, что уравнение Эйри обладает решением $v(\xi)$,
имеющим в секторе $|\arg\xi|\leqslant\pi-\delta_0$ ($\delta_0>0$ ---
произвольное фиксированное число) асимптотику
\begin{equation}
v(\xi)=\frac{e^{-\frac{2}{3}\xi^{3/2}}}{2\sqrt{\pi}\xi^{1/4}}
\left(1+\OB{\frac{1}{|\xi|^{3/2}}}\right),\quad \xi\to\infty, \label{eq_aias}
\end{equation}
где выбираются главные (положительные) ветви корней при $\xi>0$. Это решение
часто называют функцией Эйри--Фока. Нам потребуется найти более широкую
область, в которой асимптотика \eqref{eq_aias} сохраняется.

\begin{Lemma}
\label{lm_aias} Для функции Эйри--Фока $v(\xi)$ асимптотика \eqref{eq_aias}
сохраняется в области
$\Omega=\left\{\xi:\;|\arg\xi|\LE\pi-\frac{3\ln{|\xi|}}{4|\xi|^{3/2}}\right\}$.
При этом остаточный член мажорируется величиной $C|\xi|^{-3/2}$ равномерно по
$\xi\in\Omega$.
\end{Lemma}
\begin{Proof}
Справедливость формулы \eqref{eq_aias} в области $\Omega$ доказывается на
основе анализа тождества
\begin{equation}
\label{ident}
  v(\xi)= e^{-\frac{\pi i}{3}}v\left(e^{ \frac{2\pi i}{3}}\xi\right)+ e^{
  \frac{\pi i}{3}}v\left(e^{-\frac{2\pi i}{3}}\xi\right).
\end{equation}
В самом деле, формулу \eqref{eq_aias} можно использовать для вычисления
$v(e^{\pm\frac{2\pi i}{3}}\xi)$ при $\xi=-\rho e^{i\delta}$,
$\frac{\ln{|\xi|^\varkappa}}{|\xi|^{3/2}}\LE|\delta|\LE\delta_0$, $\rho>1$:
$$
  e^{\mp\frac{\pi i}{3}}v\left(e^{\pm\frac{2\pi i}{3}}\xi\right)=
  \frac{e^{-i\left(\frac{\delta}{4}\pm\frac{\pi}{4}\right)}}{2\sqrt{\pi}\rho^{1/4}}
  e^{\pm i\frac{2}{3}\rho^{3/2}e^{i\frac{3\delta}{2}}}
  \left(1+\OB{\frac{1}{\rho^{3/2}}}\right).
$$
Отсюда видно, что при $\delta>0$ первое слагаемое в тождестве \eqref{ident}
мало по сравнению со вторым и им можно пренебречь, если
$$
  e^{-i\frac{2}{3}\rho^{3/2}e^{i\frac{3\delta}{2}}}\OB{\frac{1}{\rho^{3/2}}}=
  e^{i\frac{2}{3}\rho^{3/2}e^{i\frac{3\delta}{2}}}.
$$
Это соотношение выполено при достаточно малых $\delta_0=\delta_0(\varkappa)$, так
как
$$
  \rho^{3/2}\left|e^{i\frac{4}{3}\rho^{3/2}e^{i\frac{3\delta}{2}}}\right|=
  \rho^{3/2}e^{-\rho^{3/2}\frac{4}{3}\sin{\frac{3\delta}{2}}}<
  \rho^{3/2}e^{-\rho^{3/2}\frac{3}{2}\frac{\delta}{\varkappa}}\LE
  \rho^{3/2}e^{-\frac{3}{2}\ln\rho}=
  1.
$$
Таким образом,
\begin{gather*}
  v(\xi)=e^{\frac{\pi i}{3}}
  \frac{e^{-\frac{2}{3}\left(e^{-\frac{2\pi i}{3}}\xi\right)^{3/2}}}
  {2\sqrt{\pi}\left(e^{-\frac{2\pi i}{3}}\xi\right)^{1/4}}
  \left(1+\OB{\frac{1}{|\xi|^{3/2}}}\right)=\\
  =\frac{e^{-\frac{2}{3}\xi^{3/2}}}{2\sqrt{\pi}\xi^{1/4}}
  \left(1+\OB{\frac{1}{|\xi|^{3/2}}}\right).
\end{gather*}
Случай $\delta<0$ рассматривается аналогично.
\end{Proof}

Рассмотрим область $D_\eps$ (см. рис. \ref{pic_Deps}), ограниченную слева и
справа прямыми $\RE{\lambda}=\pm 1$, а сверху --- прямыми, проходящими через
точки $\pm 1$ и точку
$d_\eps=-i\left(\frac{1}{\sqrt{3}}+\eps^{1/2}\ln\eps^{-\theta}\right)$, где
$\theta>\frac{1}{3}(3/4)^{3/4}$~--- фиксированное число. Положим
\begin{equation*}
f(\lambda)=\int\limits_{-1}^{1}e^{-i\frac{\pi}{4}}\sqrt{x-\lambda}\,dx,
\end{equation*}
где ветвь квадратного корня выбирается такой, что $f(-i/\sqrt{3})>0$.
Обозначим
$$
\Lambda_\alpha=\{\lambda\in\mathbb{C}:\; |\arg\lambda|<\alpha\}.
$$

\begin{figure}[tbp]
\unitlength=60mm
\begin{picture}(2.2,2.1)(-1.4,-2)
\input{pic_misc/Deps.tex}
\end{picture}
\caption{\small Область справедливости асимптотик в окрестности луча
$[-i/\sqrt{3},-i\infty)$} \label{pic_Deps}
\end{figure}

\begin{Lemma}
\label{lm2}
Функция $f(\lambda)$ обладает следующими свойствами:

$1^\circ.\;${}$f(\lambda)$ голоморфна в нижней полуплоскости и вещественна при
$\lambda=-i\mu,\;\mu>0$.

$2^\circ.\;${}Значения $f(\lambda)$ и $-i f'(\lambda)$ лежат в секторе
$\Lambda_{\pi/6}$, а $f''(\lambda)$ в секторе $\Lambda_{\pi/2}$ при всех
$\lambda\in D_\eps$.

$3^\circ.\;${}$f(\lambda)$ монотонно возрастает при $\lambda\to -i\infty$ по
отрицательной мнимой оси в области $D_\eps$.

$4^\circ.\;${}$f(\lambda)=2\sqrt{i\lambda}+\OB{\frac{1}{|\lambda|^{3/2}}},\;
   f'(\lambda)=\frac{i}{\sqrt{i\lambda}}+\OB{\frac{1}{|\lambda|^{3/2}}}$
   при $\lambda\in D_\eps$.

$5^\circ.\;${}$|f(\lambda_1)-f(\lambda_2)|\GE
   \gamma\frac{|\lambda_1-\lambda_2|}{\sqrt{\max\{|\lambda_1|,|\lambda_2|\}}}$
   при $\lambda_1,\lambda_2\in D_\eps$, где постоянная $\gamma>0$ не зависит от
   $\lambda_1$, $\lambda_2$ и $\eps$.

$6^\circ.\;${}$f(d_\eps)<\RE{f(\lambda)}$ при $\lambda\in D_\eps$.

$7^\circ.\;${}Функция $f(\lambda)$, $\lambda\in D_\eps$, не принимает
вещественных значений вне мнимой оси.
\end{Lemma}
\begin{VarProof}
удобнее провести в терминах переменной $\mu=i\lambda$ и функции
$g(\mu)=f(\lambda)$, которая имеет вид
$$
  g(\mu)=\int\limits_{-1}^{1}\sqrt{\mu-ix}\,dx.
$$

Свойство $1^\circ$ следует из того, что при положительных значениях $\mu$
$$
  \overline{g(\mu)}=\int\limits_{-1}^{1}\sqrt{\mu+ix}\,dx=
  \int\limits_{-1}^{1}\sqrt{\mu-i\tilde{x}}\,d\tilde{x}=g(\mu),\quad
  \tilde{x}=-x.
$$

Чтобы доказать свойство $2^\circ$, вычислим производные:
\begin{gather*}
  -if'(\lambda)=g'(\mu)=\frac{1}{2}\int\limits_{-1}^{1}\frac{1}{\sqrt{\mu-ix}}\,dx,\\
  f''(\lambda)=-g''(\mu)=
  \frac{1}{4}\int\limits_{-1}^{1}\frac{1}{(\mu-ix)^{3/2}}\,dx,
\end{gather*}
и заметим, что $\sqrt{\mu-ix}\in\Lambda_{\pi/6}$,
$\;\frac{1}{\sqrt{\mu-ix}}\in\Lambda_{\pi/6}$,
$\;\frac{1}{(\mu-ix)^{3/2}}\in\Lambda_{\pi/2}$ и интегралы от этих выражений
лежат в соответствующих секторах.

Свойство $3^\circ$ следует из $1^\circ$ и $2^\circ$, поскольку $g'(\mu)>0$ при $\mu>0$,
$-i\mu\in D_\eps$.

Для доказательства свойства $4^\circ$ вычислим интеграл:
\begin{gather*}
  g(\mu)=\left.\frac{2}{3}i(\mu-ix)^{3/2}\right|_{-1}^{1}=
  \frac{2}{3}i\mu^{3/2}\left[\left(1-\frac{i}{\mu}\right)^{3/2}-
  \left(1+\frac{i}{\mu}\right)^{3/2}\right]= \\
  =\frac{2}{3}i\mu^{3/2}\left[
  1-\frac{3}{2}\frac{i}{\mu}-\frac{3}{8}\frac{1}{\mu^2}
  -1-\frac{3}{2}\frac{i}{\mu}+\frac{3}{8}\frac{1}{\mu^2}+\OB{\frac{1}{\mu^3}}
  \right]=\\
  =2\sqrt{\mu}+\OB{\frac{1}{\mu^{3/2}}}.
\end{gather*}
Аналогично
\begin{gather*}
  g'(\mu)=\left.\lefteqn{\phantom{\frac{2}{3}}}i(\mu-ix)^{1/2}\right|_{-1}^{1}=
  i\mu^{1/2}\left[\left(1-\frac{i}{\mu}\right)^{1/2}-
  \left(1+\frac{i}{\mu}\right)^{1/2}\right]= \\
  =i\mu^{1/2}\left[
  1-\frac{1}{2}\frac{i}{\mu}+\frac{1}{8}\frac{1}{\mu^2}
  -1-\frac{1}{2}\frac{i}{\mu}-\frac{1}{8}\frac{1}{\mu^2}+\OB{\frac{1}{\mu^3}}
  \right]
  =\frac{1}{\sqrt{\mu}}+\OB{\frac{1}{\mu^{5/2}}}.
\end{gather*}
Докажем свойство $5^\circ$. Рассмотрим область $iD_0$ и докажем, что при
$\mu_1,\mu_2\in iD_0$
$$
  |g(\mu_1)-g(\mu_2)|\GE \gamma\frac{|\mu_1-\mu_2|}{\sqrt{\max\{|\mu_1|,|\mu_2|\}}}.
$$
Из свойств $2^\circ$, $4^\circ$ следует, что
$$
  \RE g'(\mu)>\frac{\gamma}{\sqrt{|\mu|}},\quad \mu\in iD_0.
$$
В самом деле, поскольку $g'(\mu)\in\Lambda_{\pi/6}$ и
$g'(\mu)\sim\frac{1}{\sqrt{\mu}}$, $\mu\to\infty$, то при достаточно больших
$|\mu|$
$$
  \RE g'(\mu)\GE \frac{\sqrt{3}}{2}|g'(\mu)|>\frac{1}{2\sqrt{|\mu|}},
$$
а при малых $|\mu|$ константу $\gamma>0$ можно подобрать, так как
$\sqrt{|\mu|}\RE g'(\mu)>0$ в замыкании области $iD_0$. Таким образом, на
отрезке $[\mu_1,\mu_2]$
$$
  \RE g'(\mu)>\frac{\gamma}{\sqrt{\max\{|\mu_1|,|\mu_2|\}}},
$$
поэтому найдется число $\xi$, для которого
$$
  g(\mu_1)-g(\mu_2)=\xi(\mu_1-\mu_2),\quad
  \RE\xi>\frac{\gamma}{\sqrt{\max\{|\mu_1|,|\mu_2|\}}}.
$$
Следовательно,
$$
  |g(\mu_1)-g(\mu_2)|=|\xi|\cdot|\mu_1-\mu_2|>
  \RE\xi\cdot|\mu_1-\mu_2|>\frac{\gamma\cdot|\mu_1-\mu_2|}{\sqrt{\max\{|\mu_1|,|\mu_2|\}}}.
$$
Докажем свойство $6^\circ$. Поскольку $g'(\mu)\in\Lambda_{\pi/6}$, то
найдется комлексное число $\xi\in\Lambda_{\pi/6}$, такое, что
$$
  g(\mu)-g(i d_\eps)=\xi(\mu-i d_\eps).
$$
И так как $\mu-i d_\eps=i(\lambda-d_\eps)\in\Lambda_{\pi/3}$, то $g(\mu)-g(i
d_\eps)\in\Lambda_{\pi/2}$, то есть
$$
  \RE(g(\mu)-g(i d_\eps))=
  \RE f(\lambda)-f(d_\eps)>
  0.
$$
Осталось доказать свойство $7^\circ$. Предположим, $g(\mu)\in\mathbb{R}$.
Тогда
$$
  g(\mu)=\overline{g(\mu)}=
  \int\limits_{-1}^{1}\sqrt{\bar{\mu}+ix}\,dx=
  \int\limits_{-1}^{1}\sqrt{\bar{\mu}-i\tilde{x}}\,d\tilde{x}=g(\bar{\mu}),\quad
  \tilde{x}=-x.
$$
Отсюда и из свойства $5^\circ$ следует, что $\mu=\bar{\mu}$, то есть $\mu$
вещественно.
\end{VarProof}

Из леммы \ref{lm2} (свойства $3^\circ$, $7^\circ$) следует, что найдется
натуральное число $k_0$ такое, что уравнение
$$
f(\lambda)=\pi k\eps^{1/2}
$$
имеет в области $D_\eps$ в точности одно чисто мнимое решение $-i\rho_k$ при всех целых
$k\GE k_0$. Очевидно,
$$
id_\eps<\rho_{k_0}<\rho_{k_0+1}<\ldots.
$$
\begin{Theorem}
\label{th_orsommod_ray} При достаточно малых $\eps>0$ все собственные
значения задачи \eqref{eq_model1}, \eqref{eq_model2}, лежащие в области
$D_\eps$, являются простыми и образуют монотонную последовательность чисел на
отрицательной мнимой оси. При этом на интервале
$\left(d_\eps,-i(\rho_{k_0}+O(\eps))\right)$ может находиться не более двух
собственных значений, а все остальные собственные значения $\lambda_k$ задачи
\eqref{eq_model1}, \eqref{eq_model2} в области $D_\eps$ имеют вид
$$
\lambda_k=-i(\rho_k+\rho_k^{-1}O(\eps)),\quad k=k_0+1,\;k_0+2,\;\ldots.
$$
\end{Theorem}
\begin{Proof}
Характеристический определитель задачи \eqref{eq_model1}, \eqref{eq_model2}
имеет вид (см.~\cite{Shkalikov_1997})
$$
v(\xi_1)v\left(e^{-\frac{2\pi i}{3}}\xi_2\right)- v\left(e^{-\frac{2\pi
i}{3}}\xi_1\right)v(\xi_2)=0.
$$
Проверим, что переменные $\; \xi_j$,\; $e^{-\frac{2\pi i}{3}}\xi_j$, $\;
j=1,2$, остаются в области, указанной в лемме \ref{lm_aias}, если $\lambda\in
D_\eps$:
\begin{gather*}
  \frac{\pi}{3}<
  \frac{\pi}{6}+\arg(1-\lambda)\LE
  \arg\xi_j \LE
  \frac{\pi}{6}+\arg(-1-\lambda)<
  \pi, \\
  -\frac{\pi}{3}<\arg e^{-\frac{2\pi i}{3}}\xi_j<\frac{\pi}{3}.
\end{gather*}
Для $\arg(-1-\lambda)$ необходимо сделать более точную оценку при малых
$\eps$:
\begin{gather*}
  \arg(-1-\lambda)<\arg(-1-d_\eps)=\\
  =\pi-\arctg\left(\frac{1}{\sqrt{3}}+\eps^{1/2}\ln\eps^{-\theta}\right)<
  \pi-\arctg\frac{1}{\sqrt{3}}-\frac{3}{4}\tilde{\theta}\eps^{1/2}\ln\eps^{-1},
\end{gather*}
где $\theta>\tilde{\theta}>\frac{1}{3}(3/4)^{3/4}$; при этом
$$
  \frac{\ln|\xi_j|^\varkappa}{|\xi_j|^{3/2}}=
  \frac{\varkappa(\frac{1}{3}\ln\eps^{-1}+\ln|1-\lambda|)}{\eps^{-1/2}|1-\lambda|^{3/2}}<
  \frac{\tilde{\varkappa}\cdot\frac{1}{3}\ln\eps^{-1}}{\eps^{-1/2}(2/\sqrt{3})^{3/2}}=
  \frac{1}{3}(3/4)^{3/4}\tilde{\varkappa}\eps^{1/2}\ln\eps^{-1},
$$
где $\tilde{\varkappa}>\varkappa>\frac{3}{4}$; таким образом,
\begin{gather*}
  \arg\xi_j\LE\frac{\pi}{6}+\arg(-1-\lambda)<
  \pi-\frac{3}{4}\tilde{\theta}\eps^{1/2}\ln\eps^{-1}<\\
  <\pi-\frac{3/4}{\tilde{\varkappa}}\frac{\tilde{\theta}}{\frac{1}{3}(3/4)^{3/4}}
  \frac{\ln|\xi_j|^\varkappa}{|\xi_j|^{3/2}}<
  \pi-\frac{\ln|\xi_j|^\varkappa}{|\xi_j|^{3/2}}
\end{gather*}
если выбрать $\varkappa$, $\tilde{\varkappa}$, $\tilde{\theta}$ так, чтобы
$$
  1<
  \frac{\varkappa}{3/4}<
  \frac{\tilde{\varkappa}}{3/4}<
  \frac{\tilde{\theta}}{\frac{1}{3}(3/4)^{3/4}}<
  \frac{\theta}{\frac{1}{3}(3/4)^{3/4}}.
$$

Воспользуемся асимптотикой леммы \ref{lm_aias} для вычисления характеристического
определителя:
\begin{align*}
  &\frac{e^{-\frac{2}{3}\left(\xi_1^{3/2}+(e^{-\frac{2\pi i}{3}}{\xi_2)}^{3/2}\right)}}
  {\xi_1^{1/4}(e^{-\frac{2\pi i}{3}}{\xi_2)}^{1/4}}
  \left(1+|\lambda|^{-3/2}\OB{\eps^{1/2}}\right)=\\
  &\qquad\qquad\qquad\qquad
  =\frac{e^{-\frac{2}{3}\left(\xi_2^{3/2}+(e^{-\frac{2\pi i}{3}}{\xi_1)}^{3/2}\right)}}
  {\xi_2^{1/4}(e^{-\frac{2\pi i}{3}}{\xi_1)}^{1/4}}
  \left(1+|\lambda|^{-3/2}\OB{\eps^{1/2}}\right).
\end{align*}
Поскольку $\frac{\pi}{3}<\arg\xi_j<\pi$, то
$$
  (e^{-\frac{2\pi i}{3}}{\xi_j)}^{3/2}=-{\xi_j}^{3/2},\qquad
  (e^{-\frac{2\pi i}{3}}{\xi_j)}^{1/4}=e^{-\frac{\pi i}{6}}{\xi_j}^{1/4},
$$
и характеристическое уравнение приводится к виду
\begin{gather}
  e^{\frac{4}{3}(\xi_2^{3/2}-\xi_1^{3/2})}=
  1+|\lambda|^{-3/2}\OB{\eps^{1/2}}, \label{eq_char_model} \\
  e^{2i\eps^{-1/2}
  \cdot\frac{2}{3}e^{-\frac{\pi i}{4}}\left[(1-\lambda)^{3/2}-(-1-\lambda)^{3/2}\right]}=
  1+|\lambda|^{-3/2}\OB{\eps^{1/2}}, \notag\\
  e^{2i\eps^{-1/2}f(\lambda)}=1+|\lambda|^{-3/2}\OB{\eps^{1/2}}.\notag
\end{gather}
После логарифмирования получаем:
\begin{equation}
\label{main_eq}
  f(\lambda)=\pi k\eps^{1/2}+|\lambda|^{-3/2}\OB{\eps},\quad
  k=0,\;\pm 1,\;\pm 2,\;\ldots.
\end{equation}
Это уравнение не может иметь решений при $k\LE k_0-2$ в силу
свойства $6^\circ$:
\begin{gather*}
  |f(\lambda)-\pi k\eps^{1/2}|\GE
  |f(\lambda)|-\pi k\eps^{1/2}\GE
  \RE f(\lambda)-\pi k\eps^{1/2}>\\
  >f(d_\eps)-\pi k\eps^{1/2}\GE
  \pi(k_0-1)\eps^{1/2}-\pi(k_0-2)\eps^{1/2}=
  \pi\eps^{1/2}.
\end{gather*}
При $k\GE k_0+1$ обозначим $B_\eps^{(k)}$ круг с центром в точке $-i\rho_k$ и
радиусом $c\rho_k^{-1}\eps$, $c>0$. В силу ограниченности $f'(\lambda)$
$$
  |\rho_k-\rho_{k-1}|>\tilde{c}\cdot\pi\eps^{1/2},\; \tilde{c}>0.
$$
поэтому круги $B_\eps^{(k)}$ не пересекаются и полностью содержатся в области
$D_\eps$ при малых $\eps>0$. Докажем, что уравнение \eqref{main_eq} имеет
ровно одно решение $\lambda_k\in D_\eps$, которое лежит внутри круга
$B_\eps^{(k)}$. Это следует из теоремы Руше и оценки для $\lambda\in
D_\eps\backslash B_\eps^{(k)}$:
$$
  |f(-i\rho_k)-f(\lambda)|\GE
  \gamma\frac{|-i\rho_k-\lambda|}{\sqrt{\max\{\rho_k,|\lambda|\}}}>
  \frac{c\gamma}{\sqrt{8}}|\lambda|^{-3/2}\eps>
  |\lambda|^{-3/2}\OB{\eps}.
$$
Первое неравенство в этой цепочке совпадает со свойством $5^\circ$, третье
выполнено при выборе достаточно большой константы $c$, а для доказательства
второго неравенства рассмотрим два случая. Если $|\lambda|\LE
\frac{1}{2}\rho_k$, то
\begin{gather*}
  \frac{|-i\rho_k-\lambda|}{\sqrt{\max\{\rho_k,|\lambda|\}}}=
  \frac{|-i\rho_k-\lambda|}{\sqrt{\rho_k}}=
  \sqrt{\rho_k}\left|1-\frac{i\lambda}{\rho_k}\right|\GE\\
  \GE\sqrt{\rho_k}\left(1-\frac{|\lambda|}{\rho_k}\right)>
  \sqrt{\frac{1}{\sqrt{3}}}\cdot\frac{1}{2}>
  c|\lambda|^{-3/2}\eps.
\end{gather*}
Если же $|\lambda|>\frac{1}{2}\rho_k$, то
$$
  \frac{|-i\rho_k-\lambda|}{\sqrt{\max\{\rho_k,|\lambda|\}}}\GE
  \frac{c\rho_k^{-1}\eps}{\sqrt{2|\lambda|}}>
  \frac{c}{\sqrt{8}}|\lambda|^{-3/2}\eps.
$$
Остается заметить, что в силу симметрии собственные значения $\lambda_k$ лежат
на мнимой оси, а случаи $k=k_0-1,k_0$ сводятся к рассмотренному, если
расширить немного область $D_\eps$ за счет уменьшения параметра $\theta$.
\end{Proof}

Асимптотика собственных значений вблизи отрезков $[\pm 1,\node]$ была
получена в \cite{Shkalikov_1997}. Здесь мы приведем уточненные формулы (ср.
\cite{Shkalikov_1997}), которые получаются при анализе известных асимптотик
функций Эйри (см. \cite{Olver_1974}) с привлечением леммы \ref{lm_aias}, как
и при доказательстве теоремы \ref{th_orsommod_ray}.

\begin{Theorem}
\label{th_orsommod_seg}
Собственные значения задачи \eqref{eq_model1}, \eqref{eq_model2} симметричны относительно мнимой оси и в
окрестности отрезка $[-1,\node]$ имеют асимптотику
\begin{align*}
&\lambda_k=(e^{-i\frac{\pi}{6}}t_k-1)
\left(1+\OB{e^{-\eps^{-1/2}\varphi(t_k)}}\right),\\
&\varphi(t)=\frac{4}{3}\RE{(2e^{i\frac{\pi}{6}}-t)^{3/2}},
\end{align*}
где $t_k$ --- нули функции Эйри--Фока $v(-\eps^{-1/3}t)=0$, когда
$t\in\left[0,\frac{2}{\sqrt{3}}-\eps^{1/2}\ln\eps^{-\theta}\right)$. При этом
\begin{equation}
\label{eq_tk}
  {(\eps^{-1/3}t_k)}^{3/2}=
  \frac{3\pi}{2}\left(k-\frac{1}{4}\right)+\OB{\frac{1}{k}},\quad k=1,2,\ldots.
\end{equation}
\end{Theorem}
\begin{Proof}
Для вывода формулы \eqref{eq_tk} в \cite{Shkalikov_1997} использовалось
следующее асимптотическое соотношение для вещественных нулей функции
Эйри--Фока $v(-x)$:
\begin{equation}
\label{eq_xk}
  x_k\sim \left[\frac{3\pi}{2}\left(k+k_0-\frac{1}{4}\right)\right]^{2/3},\quad
  k=1,2,\ldots.
\end{equation}
Это соотношение можно уточнить с помощью следующей асимптотической формулы,
известной из \cite{Olver_1974}:
$$
  \theta(x)=\frac{2}{3}|x|^{3/2}+\frac{\pi}{4}+\OB{\frac{1}{|x|^{3/2}}}, \quad
  x\to -\infty,
$$
где по определению
\begin{gather*}
  \theta(x)=\left\{
  \begin{array}{ll}
    \frac{\pi}{4}, & x\GE c, \\
    \arctg\frac{\Ai(x)}{\Bi(x)}, & x\LE c,
  \end{array}
  \right.\\
  c=\max\{x:\Ai(x)=\Bi(x)\},\quad -1<c<0.
\end{gather*}
Отсюда
\begin{gather*}
  \theta(-x_k)=\pi k,\quad k=1,2,\ldots,\\
  \frac{2}{3}x_k^{3/2}+\frac{\pi}{4}+\OB{\frac{1}{x_k^{3/2}}}=\pi k,\\
  x_k^{3/2}=\frac{3\pi}{2}\left(k-\frac{1}{4}\right)+\OB{\frac{1}{k}},\quad
  k=1,2,\ldots.
\end{gather*}
Таким образом, мы уточнили формулу \eqref{eq_xk} и вычислили $k_0=0$.
\end{Proof}

\begin{Note}
В \cite{Shkalikov_1997} вместо условия
$t_k\in\left[0,\frac{2}{\sqrt{3}}-\eps^{1/2}\ln\eps^{-\theta}\right)$ ошибочно
фигурирует условие $t_k\in\left[0,\frac{2}{\sqrt{3}}-c\eps^{1/2}\right)$.
\end{Note}

Предельное множество концентрации собственных значений
$$
T=\left[-1,\node\right) \cup \left[1,\node\right) \cup
\left(\node,-i\infty\right)
$$
названо в \cite{Shkalikov_1997} \lk спектральным галстуком\pk. Функцию
$N(\lambda,\eps)$, определенную при $\lambda\in T$, назовем функцией
распределения собственных значений задачи \eqref{eq_model1},
\eqref{eq_model2}, если для всех $\lambda_1$, $\lambda_2$, принадлежащих
любой связной компоненте множества $T$ (узел $\node$ исключается), число
собственных значений на отрезке $[\lambda_1,\lambda_2]$ (или в малом
прямоугольнике, стороны которого проходят через $\lambda_1$ и $\lambda_2$)
равно $|N(\lambda_2,\eps)-N(\lambda_1,\eps)|$.

\begin{Theorem}
\label{th_orsommod_Nlam}
Функция распределения собственных значений задачи \eqref{eq_model1},
\eqref{eq_model2} при $\eps\to 0$ на сегментах $[\pm 1,\node)$ имеет вид
$$
N(\lambda,\eps)=\pm\frac{\eps^{-1/2}}{\pi}
\int\limits_{\lambda}^{\pm 1}e^{-i\frac{\pi}{4}}\sqrt{x-\lambda}\,dx
+\OB{1},
$$
а на луче $(\node,-i\infty)$
$$
N(\lambda,\eps)=\frac{\eps^{-1/2}}{\pi}
\int\limits_{-1}^{1}e^{-i\frac{\pi}{4}}\sqrt{x-\lambda}\,dx
+\OB{1}.
$$
\end{Theorem}
\begin{Proof}
По теореме \ref{th_orsommod_seg} на отрезке $[-1,\node]$ имеем:
$$
  {(\eps^{-1/3}t_k)}^{3/2}=\frac{3\pi k}{2}+O(1).
$$
Отсюда
$$
  k=
  \frac{\eps^{-1/2}}{\pi}\cdot\frac{2}{3}t_k^{3/2}+O(1)=
  \frac{\eps^{-1/2}}{\pi}\cdot\frac{2}{3}|-1-\lambda_k|^{3/2}+O(1).
$$
С другой стороны, при $\lambda\in[-1,\node]$
\begin{gather*}
  -\int\limits_{\lambda}^{-1}e^{-i\frac{\pi}{4}}\sqrt{x-\lambda}\,dx=
  -\left.\frac{2}{3}e^{-i\frac{\pi}{4}}(x-\lambda)^{3/2}\right|_\lambda^{-1}=\\
  \frac{2}{3}e^{i\frac{3\pi}{4}}(-1-\lambda)^{3/2}=
  \frac{2}{3}|-1-\lambda|^{3/2}.
\end{gather*}
Вычислим теперь функцию распределения на луче $(\node,-i\infty)$,
которую можно в данном случае определить как количество собственных значений,
мнимая часть которых больше $\IM\lambda$. Пусть
$f(\lambda)=\pi\eps_0^{1/2}(k_0\pm \frac{1}{2})$, тогда
\begin{gather*}
  N(\lambda,\eps)-N(\lambda,\eps_0)=k-k_0\pm 2,\qquad\mbox{где}\quad
  f(\lambda)=\pi\eps^{1/2}(k\pm\frac{1}{2}), \\
  N(\lambda,\eps)=\frac{\eps^{-1/2}}{\pi}f(\lambda)+O(1).
\end{gather*}
\end{Proof}

\begin{Note}
Теорема \ref{th_orsommod_Nlam} была получена авторами в 1998 году, а
впоследствии обобщена на случай, когда вместо $q(x)=x$ участвует функция
$q(x)=x^2$ или аналитическая монотонная функция (см. Туманов, Шкаликов
\cite{Tumanov_Shkalikov_2002}, Шкаликов \cite{Shkalikov_2001}).
\end{Note}

\begin{Theorem}
\label{th_orsommod_Ndel} Количество собственных значений задачи
\eqref{eq_model1}, \eqref{eq_model2} в круге малого радиуса $\delta$ с
центром в точке-узле $\node$ равно (при $\eps\to 0$)
$$
N_\delta(\eps)=\frac{\eps^{-1/2}}{\pi} \left(
f\left(\node-i\delta\right)-\frac{4}{3}\left(\frac{2}{\sqrt{3}}-\delta\right)^{3/2}
\right)+\OB{1}.
$$
В круге уменьшающегося радиуса $\eps^{1/2}\ln\eps^{-\theta}$, при
фиксированном $\theta>\frac{1}{3}(3/4)^{3/4}$, число собственных значений
равно
$$
N_{\eps^{1/2}\ln\eps^{-\theta}}(\eps)=
\frac{2^{1/2}3^{3/4}}{\pi}\ln\eps^{-\theta}+\OB{1},\quad \eps\to 0.
$$
\end{Theorem}
\begin{Proof}
  Анализ доказательства теоремы \ref{th_orsommod_Nlam} дает формулу для $N_\delta(\eps)$ при
  фиксированном $\delta>0$. Устремляя $\delta$ к нулю, получим:
  \begin{gather*}
  f\left(\node-i\delta\right)=
  f\left(\node\right)-
  i\delta f'\left(\node\right)+
  O(\delta^2), \\
  f\left(\node\right)=
  \left.\frac{2}{3}e^{-i\frac{\pi}{4}}\left(x+\frac{i}{\sqrt{3}}\right)^{3/2}\right|_{-1}^1=
  2\cdot\frac{2}{3}\left(\frac{2}{\sqrt{3}}\right)^{3/2}=
  \frac{4}{3}\left(\frac{2}{\sqrt{3}}\right)^{3/2}, \\
  f'\left(\node\right)=
  \left.-e^{-i\frac{\pi}{4}}\left(x+\frac{i}{\sqrt{3}}\right)^{1/2}\right|_{-1}^1=
  -\left.\left(\frac{1}{\sqrt{3}}-ix\right)^{1/2}\right|_{-1}^1=\\
  =\left.\left(\frac{1}{\sqrt{3}}+ix\right)^{1/2}\right|_{-1}^1=
  i\left(\frac{2}{\sqrt{3}}\right)^{1/2}, \\
  \frac{4}{3}\left(\frac{2}{\sqrt{3}}-\delta\right)^{3/2}=
  \frac{4}{3}\left(\frac{2}{\sqrt{3}}\right)^{3/2}
  -2\delta\left(\frac{2}{\sqrt{3}}\right)^{1/2}+
  O(\delta^2), \\
  N_\delta(\eps)=
  \frac{\eps^{-1/2}}{\pi} \left(
  3\delta\left(\frac{2}{\sqrt{3}}\right)^{1/2}+O(\delta^2)
  \right)+\OB{1}=\\
  =\frac{2^{1/2}3^{3/4}}{\pi}\delta\eps^{-1/2}+\OB{\delta^2\eps^{-1/2}}+\OB{1}
  \end{gather*}
  Теорема \ref{th_orsommod_ray} позволяет взять $\delta=\eps^{1/2}\ln\eps^{-\theta}$:
$$
  N_{\eps^{1/2}\ln\eps^{-\theta}}(\eps)=
  \frac{2^{1/2}3^{3/4}}{\pi}\ln\eps^{-\theta}+\OB{1}.
$$
\end{Proof}

\section{Решение уравнения Орра--Зоммерфельда}
\label{par32}

Перейдем теперь к изучению самой задачи Орра--Зоммерфельда с линейным профилем
\eqref{orsom1}, \eqref{orsom2}. Перепишем ее в виде:
\begin{gather}
-i\eps w''=(x-\lambda)w, \label{model0}\\
 y''-\alpha^2y=w(x), \label{eq_w}\\
 y(-1)=y'(-1)=0, \label{model0_b1}\\
 y(1)=y'(1)=0. \label{model0_b2}
\end{gather}
Функцию Грина для задачи \eqref{eq_w}, \eqref{model0_b1} найдем методом
вариации постоянных:
\begin{align*}
y(x)&=c_1(x)e^{\alpha x}+c_2(x)e^{-\alpha x},\\
y'(x)&=\alpha c_1(x)e^{\alpha x}-\alpha c_2(x)e^{-\alpha x}+c'_1(x)e^{\alpha
x}+c'_2(x)e^{-\alpha x}.
\end{align*}
Получаем первое уравнение для $c'_1(x)$, $c'_2(x)$:
\begin{equation}
c'_1(x)e^{\alpha x}+c'_2(x)e^{-\alpha x}=0. \label{cslu1}
\end{equation}
Найдем вторую производную и подставим в \eqref{eq_w}:
\begin{align*}
y''(x)&=\alpha^2 c_1(x)e^{\alpha x}+\alpha^2 c_2(x)e^{-\alpha x}+
\alpha c'_1(x)e^{\alpha x}-\alpha c'_2(x)e^{-\alpha x}=\\
&=\alpha^2 c_1(x)e^{\alpha x}+\alpha^2 c_2(x)e^{-\alpha x}+w(x).
\end{align*}
Получаем второе уравнение для $c'_1(x)$, $c'_2(x)$:
\begin{equation}
c'_1(x)e^{\alpha x}-c'_2(x)e^{-\alpha x}=\frac{1}{\alpha}w(x). \label{cslu2}
\end{equation}
Решение системы линейных уравнений \eqref{cslu1}, \eqref{cslu2} дает:
\begin{gather*}
c'_1(x)=\frac{1}{\alpha}\frac{e^{-\alpha x}}{2}w(x),\quad
c'_2(x)=-\frac{1}{\alpha}\frac{e^{\alpha x}}{2}w(x),\\
y(x)=\frac{1}{\alpha}\int\limits_{x_1}^x \frac{e^{\alpha(x-t)}}{2}w(t)\,dt-
\frac{1}{\alpha}\int\limits_{x_2}^x \frac{e^{-\alpha(x-t)}}{2}w(t)\,dt.
\end{gather*}
С учетом условий \eqref{model0_b1} получаем:
\begin{equation*}
y(x)=\frac{1}{\alpha}\int\limits_{-1}^x \sh\left[\alpha(x-t)\right]w(t)\,dt,
\qquad y'(x)=\int\limits_{-1}^x \ch\left[\alpha(x-t)\right]w(t)\,dt.
\end{equation*}

Таким образом, исходную задачу \eqref{orsom1}, \eqref{orsom2} можно записать в
виде уравнения второго порядка с нелокальными краевыми условиями:
\begin{gather*}
-i\eps w''=(x-\lambda)w,\\
\int\limits_{-1}^1 \sh\left[\alpha(1-t)\right]w(t)\,dt=\int\limits_{-1}^1
\ch\left[\alpha(1-t)\right]w(t)\,dt=0.
\end{gather*}
Такая запись краевой задачи Орра-Зоммерфельда известна давно.
Характеристическое уравнение теперь запишется так:
\begin{equation*}
\begin{vmatrix}
  \int\limits_{-1}^1 \sh\left[\alpha(1-t)\right]w_1(t)\,dt &
  \int\limits_{-1}^1 \sh\left[\alpha(1-t)\right]w_2(t)\,dt \\
  \int\limits_{-1}^1 \ch\left[\alpha(1-t)\right]w_1(t)\,dt &
  \int\limits_{-1}^1 \ch\left[\alpha(1-t)\right]w_2(t)\,dt
\end{vmatrix}
=0,
\end{equation*}
где $w_1(x)$, $w_2(x)$ --- два независимых решения уравнения \eqref{model0}.
Замена
$$
  \xi_j=\omega_j\,\sigma^{-1}(x-\lambda),\qquad j=0,1,2,
$$
где
\begin{align*}
  &\omega_j=e^{i(\frac{2\pi}{3}j+\frac{\pi}{6})},\\
  &\sigma=\eps^{1/3},
\end{align*}
приводит \eqref{model0} к модельному уравнению Эйри
$$
  w''(\xi)=\xi w(\xi),
$$
поэтому можно взять
\begin{equation*}
w_n(x)=v(\xi_{j_n})=v(\omega_{j_n}\sigma^{-1}(x-\lambda)),
\end{equation*}
где $v(\xi)$ --- функция Эйри--Фока.

Производя соответствующую замену в интегралах, входящих в характеристический
определитель, можно переписать их в виде:
\begin{align*}
 &\int\limits_{-1}^1 \sh\left[\alpha(1-t)\right]w_n(t)\,dt=
  \omega_{j_n}^{-1}\sigma\int\limits_{\xi_{j_n}^-}^{\xi_{j_n}^+}g_{j_n}(\sigma z)v(z)\,dz, \\
 &\int\limits_{-1}^1 \ch\left[\alpha(1-t)\right]w_n(t)\,dt=
  \omega_{j_n}^{-1}\sigma\int\limits_{\xi_{j_n}^-}^{\xi_{j_n}^+}g_{j_n}^*(\sigma z)v(z)\,dz,
\end{align*}
где
\begin{align*}
  &\xi_j^\pm=\omega_j\,\sigma^{-1}(\pm 1-\lambda),\\
  &g_j(z)=\sh\left[\alpha(1-\lambda-\omega_j^{-1}z)\right],\\
  &g_j^*(z)=\ch\left[\alpha(1-\lambda-\omega_j^{-1}z)\right].
\end{align*}
Обозначим
\begin{align}\label{ujxij}
 u_j(\xi_j)=\int\limits_{\xi_j^0}^{\xi_j}g_j(\sigma z)v(z)\,dz, \qquad
 u_j^*(\xi_j)=\int\limits_{\xi_j^0}^{\xi_j}g_j^*(\sigma z)v(z)\,dz,
\end{align}
где $\xi_j^0$ --- некоторые числа, и окончательно запишем характеристическое
уравнение так:
\begin{equation}
\begin{vmatrix}\label{detuxi}
u_{j_1}(\xi_{j_1}^+)-u_{j_1}(\xi_{j_1}^-) &
u_{j_2}(\xi_{j_2}^+)-u_{j_2}(\xi_{j_2}^-) \\
u_{j_1}^*(\xi_{j_1}^+)-u_{j_1}^*(\xi_{j_1}^-) &
u_{j_2}^*(\xi_{j_2}^+)-u_{j_2}^*(\xi_{j_2}^-)
\end{vmatrix}
=0.
\end{equation}

Собственные значения $\lambda$ задачи Орра--Зоммерфельда, как и собственные
значения модельной задачи, симметричны относительно мнимой оси и при $\eps\to
0$ концентрируются в окрестности отрезков $[1,\node]$, $[-1,\node]$ и луча
$[\node,-i\infty)$. В следующих двух параграфах вычисляется асимптотика
собственных значений в окрестности отрезка $[-1,\node]$ и луча
$[\node,-i\infty)$.

\section{Поведение собственных значений в окрестности отрезка $[-1,\node]$}
\label{par33}

Выберем на отрезке $[-1,\node]$ точки $d_1$ и $d_2$, которые при $\eps\to 0$
приближаются к концам отрезка $-1$ и $\node$ так, что
\begin{gather}
  \eps^{-1/3}|-1-d_1|\to +\infty, \label{orsom_cond1}\\
  \left|\node-d_2\right|=2\eps^{1/2}\ln{\eps^{-\theta}},\quad
  \theta>\frac{1}{3}\left(\frac{3}{4}\right)^{3/4}. \label{orsom_cond2}
\end{gather}
Рассмотрим область $\Omega_\eps=\Omega_\eps^1\cup\Omega_\eps^2$ (см. рис.
\ref{pic_Omegaeps}), являющуюся частью полуполосы $-1\LE\RE{\lambda}\LE 0$,
$\IM{\lambda}\LE 0$, которая отделена от точки $-1$ дугой окружности с
центром $-1$, проходящей через точку $d_1$, и ограничена снизу прямой,
проходящей через точки $d_2$ и $1$.

\begin{figure}[tbp]
\unitlength=60mm
\begin{picture}(2.2,1.1)(-1.4,-1)
\input{pic_misc/Omegaeps.tex}
\end{picture}
\caption{\small Области справедливости асимптотик в окрестности отрезка
$[-1,-i/\sqrt{3}]$} \label{pic_Omegaeps}
\end{figure}

Если $\lambda$ лежит в этой области, то при $\sigma\to 0$
\begin{gather}
  \xi_j^\pm\to\infty,\qquad \sigma\xi_j^\pm=\OB{1} \label{orsom_cond1a}\\
  |\arg\xi_j^+|\LE\pi-2\frac{\ln{|\xi_j^+|^\varkappa}}{|\xi_j^+|^{3/2}},\quad
  \varkappa>\frac{3}{4}.
  \label{orsom_cond2a}
\end{gather}

Отрезок $[-1,\node]$ разбивает $\Omega_\eps$ на две части: верхнюю
$\Omega_\eps^1$ и нижнюю $\Omega_\eps^2$. В характеристическом определителе
\eqref{detuxi} положим $j_1=0$, $j_2=l$, если $\lambda\in\Omega_\eps^l$.
Число $\xi_j^0$ в формуле \eqref{ujxij} положим равным нулю:
\begin{align}\label{main_int}
 u_j(\xi_j)=\int\limits_{0}^{\xi_j}g_j(\sigma z)v(z)\,dz, \qquad
 u_j^*(\xi_j)=\int\limits_{0}^{\xi_j}g_j^*(\sigma z)v(z)\,dz,
\end{align}

\begin{Lemma}
\label{lm_uxi_prop}
 Интегралы \eqref{main_int} обладают следующими свойствами
\begin{align}
& 1^\circ.\quad\sum\limits_{j=0}^{2}u_j(\xi_j^\pm)=0. \label{prop1}\\
& 2^\circ.\quad  u(\xi)=\frac{1}{3}\left(g(0)+\OB{\xi^{-3/4}}\right)-\notag\\
& \quad\qquad\qquad\qquad\qquad
  -\frac{e^{-\frac{2}{3}\xi^{3/2}}}{2\sqrt{\pi}\xi^{3/4}}
  \left(g(\sigma\xi)+\OB{\xi^{-3/2}}\right),\quad
  \sigma\to 0.
  \label{prop2}
\end{align}
Здесь $\xi=\xi_j^+$, $\xi=\xi_1^-$ или $\xi=\xi_2^-$, а в качестве $u$, $g$
можно взять $u_j$, $g_j$ или $u_j^*$, $g_j^*$.
\end{Lemma}
\begin{Proof}
Первое свойство легко выводится из тождества:
$$
\sum\limits_{j}e^{\frac{2\pi i}{3}j}v(e^{\frac{2\pi i}{3}j}\xi)=0.
$$
В самом деле,
\begin{gather*}
  \sum\limits_{j}u_j(\xi_j^\pm)=
  \sum\limits_{j}\int\limits_{0}^{\xi_j^\pm}g_j(\sigma z)v(z)\,dz=
  \sum\limits_{j}\int\limits_{0}^{\omega_j^{-1}\xi_j^\pm}
  g_j(\sigma\omega_j w)v(\omega_j w)\omega_j\,dw=\\
  =\int\limits_{0}^{\sigma^{-1}(\pm 1-\lambda)}
  \sh\left[\alpha(1-\lambda-\sigma w)\right]\sum\limits_{j}\omega_j v(\omega_j w)\,dw=0.
\end{gather*}

Чтобы вычислить интеграл $u(\xi)$, разобьем его на две части:
$$
  u(\xi)=\int\limits_{0}^{|\xi|}g(\sigma z)v(z)\,dz+
  \int\limits_{|\xi|}^{\xi}g(\sigma z)v(z)\,dz,
$$
причем первый из этих интегралов берется по отрезку положительной полуоси, а
второй~--- по дуге окружности, не пересекающей отрицательную полуось.

Вычислим первый интеграл. Для этого заметим, что существует ограниченная
вместе со своей производной функция $h(z)$, такая, что
\begin{gather*}
  \int\limits_{0}^{|\xi|}g(\sigma z)v(z)\,dz=
  \int\limits_{0}^{|\xi|}(g(0)+\sigma z h(z))v(z)\,dz= \\
  =g(0)\int\limits_{0}^{\infty}v(z)\,dz-
  g(0)\int\limits_{|\xi|}^{\infty}v(z)\,dz+
  \int\limits_{0}^{|\xi|}\sigma z h(z)v(z)\,dz.
\end{gather*}
Интеграл от функции Эйри по положительной полуоси можно вычислить (см.
\cite[гл.~11, \S 12.2]{Olver_1974}):
\begin{equation*}
  \int\limits_{0}^{\infty}v(z)\,dz=\frac{1}{3}.
\end{equation*}
Оценим остаток:
\begin{gather*}
  g(0)\int\limits_{|\xi|}^{\infty}v(z)\,dz=
  \OB{\int\limits_{|\xi|}^{\infty}\frac{e^{-\frac{2}{3}z^{3/2}}}{z^{1/4}}\,dz}=\\
  =\OB{\int\limits_{|\xi|}^{\infty}\frac{d e^{-\frac{2}{3}z^{3/2}}}{z^{3/4}}}=
  \OB{\frac{e^{-\frac{2}{3}|\xi|^{3/2}}}{|\xi|^{3/4}}}=
  \OB{\xi^{-3/4}}; \\
  \int\limits_{0}^{|\xi|}\sigma z h(z)v(z)\,dz=
  \sigma\int\limits_{0}^{|\xi|}h(z)v''(z)\,dz= \\
  =\sigma\left( \left. h(z)v'(z) \lefteqn{\phantom\int}\right|_0^{|\xi|}-
  \int\limits_{0}^{|\xi|}h'(z)v'(z)\,dz \right)=
  \sigma\OB{1}=\OB{\xi^{-1}}.
\end{gather*}

Для вычисления второго интеграла сделаем замену
$$
  w=\frac{2}{3}(z^{3/2}-\xi^{3/2}),\qquad
  \eta=\frac{2}{3}(|\xi|^{3/2}-\xi^{3/2}).
$$
Поскольку $|z|=|\xi|$, то, с учетом условий \eqref{orsom_cond1a},
\eqref{orsom_cond2a}, можно воспользоваться леммой \ref{lm_aias}:
\begin{gather*}
  v(z)=\frac{e^{-\frac{2}{3}z^{3/2}}}{2\sqrt{\pi}z^{1/4}}
  \left(1+\OB{\frac{1}{z^{3/2}}}\right)=
  \frac{e^{-\frac{2}{3}\xi^{3/2}}}{2\sqrt{\pi}}z^{-1/4}\left(1+\OB{\xi^{-3/2}}\right)e^{-w},
  \\
  \int\limits_{|\xi|}^{\xi}g(\sigma z)v(z)\,dz=
  -\frac{e^{-\frac{2}{3}\xi^{3/2}}}{2\sqrt{\pi}}\int\limits_0^\eta
  z^{-3/4}g(\sigma z)\left(1+\OB{\xi^{-3/2}}\right)e^{-w}\,dw.
\end{gather*}
Заметим, что
\begin{gather*}
  \frac{d}{dw}\left(z^{-3/4}g(\sigma z)\right)=
  \frac{d}{dz}\left(z^{-3/4}g(\sigma z)\right)\frac{dz}{dw}=\\
  =\left(z^{-3/4}\sigma g'(\sigma z)-\frac{3}{4}z^{-7/4}g(\sigma z)\right)
  z^{-1/2}=\OB{\xi^{-9/4}},
\end{gather*}
поэтому
\begin{gather*}
  z^{-3/4}g(\sigma z)=\xi^{-3/4}g(\sigma \xi)+w\OB{\xi^{-9/4}}, \\
  \int\limits_{|\xi|}^{\xi}g(\sigma z)v(z)\,dz=
  -\frac{e^{-\frac{2}{3}\xi^{3/2}}}{2\sqrt{\pi}\xi^{3/4}}\int\limits_0^\eta
  \left(g(\sigma \xi)+w\OB{\xi^{-3/2}}\right)\left(1+\OB{\xi^{-3/2}}\right)e^{-w}\,dw.
\end{gather*}
Если $|\arg\xi|>\frac{\pi}{3}$, то
$\quad\RE\eta\GE\frac{2}{3}|\xi|^{3/2}\to+\infty,\quad
\frac{|\eta|}{\RE\eta}\LE\frac{4}{3}\frac{|\xi|^{3/2}}{\RE\eta}\LE
  2=\OB{1}$,
\begin{gather*}
  \int\limits_0^\eta g(\sigma\xi) e^{-w}\,dw=g(\sigma\xi)(1-e^{-\eta})=
  g(\sigma\xi)+\OB{\xi^{-\infty}}, \\
  \int\limits_0^\eta\left(g(\sigma\xi)O(\xi^{-3/2})+w
  O(\xi^{-3/2})\right)e^{-w}\,dw=\\
  =\int\limits_0^{\RE\eta}\xi^{-3/2}(\OB{1}+\OB{\tilde w})
  e^{-\frac{\eta}{\RE\eta}\tilde w}\,d\tilde w= \\
  =\OB{|\xi|^{-3/2}\int\limits_0^{+\infty}(1+\tilde w)e^{-\tilde w}\,d\tilde
  w}=\OB{\xi^{-3/2}}, \\
  \int\limits_{|\xi|}^{\xi}g(\sigma z)v(z)\,dz=
  -\frac{e^{-\frac{2}{3}\xi^{3/2}}}{2\sqrt{\pi}\xi^{3/4}}
  \left(g(\sigma \xi)+\OB{\xi^{-3/2}}\right).
\end{gather*}
Если $|\arg\xi|\LE\frac{\pi}{3}$, то последняя формула не имеет места, но
поскольку $\quad\RE\xi^{3/2}\GE 0,\quad e^{-\frac{2}{3}\xi^{3/2}}=\OB{1}$, то
левая и правая ее части равны $\OB{\xi^{-3/4}}$ и это есть остаточный член
доказываемой асимптотики \eqref{prop2}.
\end{Proof}

\begin{Theorem}\label{th_orsom_main}
Введем для удобства прямоугольную систему координат $(t,\gamma)$ на
комплексной плоскости $\lambda$. Начало координат расположим в точке
$\lambda=-1$, а ось $t$ направим вдоль отрезка $[-1,\node]$. Тогда
асимптотически, при $\eps\to 0$, собственные значения задачи \eqref{orsom1},
\eqref{orsom2} из области $\Omega_\eps$ располагаются на двух линиях
\begin{gather*}
  \gamma^\pm=\pm\frac{\eps^{1/2}}{t^{1/2}}\left(\ln\frac{c(\lambda)t^{3/4}}{\eps^{1/4}}+
  \OB{\frac{\eps^{1/4}}{t^{3/4}}}\right),\quad
  T_1(\eps) \LE t \LE \frac{2}{\sqrt{3}}-T_2(\eps),\\
  \mbox{где}\quad \eps^{-1/3}T_1(\eps)\to +\infty,\quad
  T_2(\eps)=2\eps^{1/2}\ln{\eps^{-\theta}},\quad
  \theta>\frac{1}{3}\left(\frac{3}{4}\right)^{3/4}.
\end{gather*}
При этом для самих собственных значений $\lambda_k=(t_k,\gamma_k)$ в
координатной плоскости $(t,\gamma)$ имеем:
\begin{equation}\label{eq_t_k}
  \frac{\left(t_k^\pm\right)^{3/2}}{\eps^{1/2}}=3\pi(k-\frac{1}{8}\mp k_0(\lambda))+
  \OB{\frac{1}{\sqrt{k}}},\qquad \gamma_k=\gamma(t_k).
\end{equation}
 Величины $c(\lambda)$ и $k_0(\lambda)$ вычисляются явно:
$$
  c(\lambda)=2\sqrt{\pi}\frac{\left|\sh\left[\alpha(1-\lambda)\right]\right|}{\sh
  2\alpha},\qquad
  k_0(\lambda)=\frac{1}{2\pi}\arg\sh\left[\alpha(1-\lambda)\right].
$$
В частности, можно утверждать, что с точностью до $\ln\eps^{-1/4}$
собственные значения в области $D_\eps$ находятся на кривых
$\gamma^\pm(t)=\pm\frac{1}{4}\eps^{1/2}t^{-1/2}\ln\eps^{-1}t^3$ в
координатной плоскости $(t,\gamma)$, причем их координаты $(t_k,\gamma_k)$
вычисляются по формулам \eqref{eq_t_k}.

Функция распределения собственных значений для каждой из кривых $\gamma^+$,
$\gamma^-$ равна половине соответствующей функции распределения для модельной
задачи на отрезке $[-1,\node]$:
$$
N^\pm(\lambda,\eps)=-\frac{\eps^{-1/2}}{2\pi}
\int\limits_{\lambda}^{-1}e^{-i\frac{\pi}{4}}\sqrt{x-\lambda}\,dx +\OB{1}.
$$
\end{Theorem}
\begin{Proof}
При $\lambda\in\Omega_\eps^l$, $l=1,2$, перепишем характеристическое
уравнение, пользуясь свойством \eqref{prop1} из леммы \ref{lm_uxi_prop}:
$$
\begin{vmatrix}
  u_0(\xi_0^+)+u_1(\xi_1^-)+u_2(\xi_2^-) &
  u_l(\xi_l^+)-u_l(\xi_l^-) \\
  u_0^*(\xi_0^+)+u_1^*(\xi_1^-)+u_2^*(\xi_2^-) &
  u_l^*(\xi_l^+)-u_l^*(\xi_l^-)
\end{vmatrix}
=0.
$$
Воспользуемся свойством \eqref{prop2} для вычисления
$u_l(\xi_l^+)-u_l(\xi_l^-)$. Поскольку
$-\frac{\pi}{3}\LE\arg\xi_l^-\LE\frac{\pi}{3}$, то функции $u_l(\xi_l^-)$
ограничены. Первое слагаемое в разложении \eqref{prop2} для $u_l(\xi_l^+)$
также ограничено, а второе слагаемое доминирует в этом разложении, так как
$$
  1=\OB{\frac{e^{-\frac{2}{3}(\xi_l^+)^{3/2}}}{(\xi_l^+)^{3/2}}}.
$$
Это следует из того, что $\frac{5\pi}{6}\LE\arg\xi_1^+\LE
\pi-\frac{3}{2}\frac{\ln{|\xi_1^+|}}{|\xi_1^+|^{3/2}},\quad
(\xi_1^+)^{3/2}=(\xi_2^+)^{3/2}$. Те же рассуждения применимы к
$u_l^*(\xi_l^+)-u_l^*(\xi_l^-)$ и мы получаем следующую оценку:
$$
  \frac{u_l(\xi_l^+)-u_l(\xi_l^-)}{u_l^*(\xi_l^+)-u_l^*(\xi_l^-)}=
  -\frac{g_l(\sigma\xi_l^+)+\OB{(\xi_l^+)^{-3/4}}}
  {g_l^*(\sigma\xi_l^+)+\OB{(\xi_l^+)^{-3/4}}}=
  -\frac{\sh 0+\OB{\sigma^{3/4}}}
  {\ch 0+\OB{\sigma^{3/4}}}=\OB{\sigma^{3/4}}.
$$
Таким образом, если $\lambda\in\Omega_\eps$, то характеристическое уравнение
можно записать так:
\begin{equation}
  u_0(\xi_0^+)+u_1(\xi_1^-)+u_2(\xi_2^-)=
  (u_0^*(\xi_0^+)+u_1^*(\xi_1^-)+u_2^*(\xi_2^-))O(\sigma^{3/4}).
  \label{char_up}
\end{equation}
  Воспользуемся еще раз леммой \ref{lm_uxi_prop}. Заметим, что
  $\xi_0^+\in[-\frac{\pi}{3},\frac{\pi}{3}]$ и
  $\xi_{3-l}^-\in[-\frac{\pi}{3},\frac{\pi}{3}]$, когда
  $\lambda\in\Omega_\eps^{3-l}$, $l=1,2$, поэтому в соответствующих асимптотических
  разложениях \eqref{prop2} можно отбросить второе слагаемое, так как оно
  мало по сравнению с первым:
\begin{align*}
  u_0(\xi_0^+)&=\frac{1}{3}(g_0(0)+O((\xi_0^+)^{-3/4})), \\
  u_{3-l}(\xi_{3-l}^-)&=\frac{1}{3}(g_{3-l}(0)+O((\xi_{3-l}^-)^{-3/4})), \\
  u_l(\xi_l^-)&=\frac{1}{3}(g_l(0)+O((\xi_l^-)^{-3/4}))-
  \frac{e^{-\frac{2}{3}(\xi_l^-)^{3/2}}}{2\sqrt{\pi}(\xi_l^-)^{3/4}}
  (g_l(\sigma\xi_l^-)+O((\xi_l^-)^{-3/2})).
\end{align*}
  Представим $\lambda$ в виде
$$
  \lambda=-1+e^{-i\frac{\pi}{6}}\eps^{1/3}\mu^{2/3}e^{i\delta},\quad
  \mu\to+\infty.
$$
Тогда
\begin{gather*}
  |\xi_j^\pm|^{-1}=\left|\omega_j\,\sigma^{-1}(\pm 1-\lambda)\right|^{-1}=
  \frac{\eps^{1/3}}{|\pm 1-\lambda|}=
  \frac{|-1-\lambda|}{|\pm 1-\lambda|}\mu^{-2/3}=O(\mu^{-2/3})
\end{gather*}
и поскольку
\begin{equation*}
  g_j(0)=\sh\left[\alpha(1-\lambda)\right],\quad g_j(\sigma\xi_j^-)=\sh
  2\alpha,
\end{equation*}
то
\begin{align*}
  u_0(\xi_0^+)&=\frac{1}{3}(\sh\left[\alpha(1-\lambda)\right]+O(\mu^{-1/2})), \\
  u_{3-l}(\xi_{3-l}^-)&=\frac{1}{3}(\sh\left[\alpha(1-\lambda)\right]+O(\mu^{-1/2})), \\
  u_l(\xi_l^-)&=\frac{1}{3}(\sh\left[\alpha(1-\lambda)\right]+O(\mu^{-1/2}))-
  \frac{e^{-\frac{2}{3}(\xi_l^-)^{3/2}}}{2\sqrt{\pi}(\xi_l^-)^{3/4}}
  (\sh 2\alpha+O(\mu^{-1})).
\end{align*}
Аналогичные разложения имеют место для $u_0^*(\xi_0^+)$,
$u_{3-l}^*(\xi_{3-l}^-)$, $u_l^*(\xi_l^-)$, поэтому характеристическое
уравнеие \eqref{char_up} сводится к следующему:
$$
  \frac{e^{-\frac{2}{3}(\xi_l^-)^{3/2}}}{2\sqrt{\pi}(\xi_l^-)^{3/4}}
  (\sh 2\alpha+O(\mu^{-1/2}))=
  \sh\left[\alpha(1-\lambda)\right]+O(\mu^{-1/2}).
$$
Далее,
$$
  \xi_l^-=\omega_l\,\sigma^{-1}(-1-\lambda)=
  -e^{-i\frac{\pi}{6}}\omega_l\mu^{2/3}e^{i\delta}=
  -e^{i\frac{2\pi l}{3}}\mu^{2/3}e^{i\delta}=
  e^{\pm i\frac{\pi}{3}}\mu^{2/3}e^{i\delta},
$$
и мы приходим к такому уравнению:
$$
  e^{\mp\frac{2}{3}i\mu e^{i\frac{3\delta}{2}}}=
  2\sqrt{\pi}\frac{\sh\left[\alpha(1-\lambda)\right]}{\sh 2\alpha}
  e^{\pm i\frac{\pi}{4}}\sqrt\mu e^{i\frac{3\delta}{4}}
  (1+O(\mu^{-1/2})).
$$
После логарифмирования получим:
\begin{align*}
  &\mp\frac{2}{3}i\mu\cos\frac{3\delta}{2}\pm\frac{2}{3}\mu\sin\frac{3\delta}{2}=\\
  &\qquad\qquad\qquad
  =\mp 2 i \pi k
  +\ln\left(c(\lambda)\sqrt{\mu}\right)
  +2 i \pi k_0(\lambda)
  \pm i\frac{\pi}{4}
  +i\frac{3\delta}{4}
  +O(\mu^{-1/2}).
\end{align*}
Равенство действительных частей дает
\begin{align}
  \frac{2}{3}\sin\frac{3\delta}{2}
  &=\pm\frac{\ln\left(c(\lambda)\sqrt{\mu}\right)}{\mu}+
  \OB{\frac{1}{\mu^{3/2}}},\notag\\
  \delta
  &=\pm\frac{\ln\left(c(\lambda)\sqrt{\mu}\right)}{\mu}+
  \OB{\frac{1}{\mu^{3/2}}}.\label{curve_delta}
\end{align}
Отсюда, в частности, следует, что $\delta=O(\mu^{-1/2}),\;
\cos\frac{3\delta}{2}=1+O(\mu^{-3/2})$, и равенство мнимых частей дает
\begin{align}
  \mp\frac{2}{3}\mu&=\mp 2\pi k+2\pi k_0(\lambda)\pm\frac{\pi}{4}
  +O(\mu^{-1/2}),\notag\\
  \mu&=3\pi(k-\frac{1}{8}\mp k_0(\lambda))
  +\OB{\frac{1}{\mu^{1/2}}}.\label{eigen_mu}
\end{align}
Мы получили заявленные результаты в полярных координатах $(\mu,\delta)$.
Остается перейти к координатам $(t,\gamma)$ по формулам
\begin{align*}
t &= \eps^{1/3}\mu^{2/3}\cos\delta,\\
\gamma &= \eps^{1/3}\mu^{2/3}\sin\delta.
\end{align*}
Обратим эти формулы:
\begin{align*}
\delta &= \gamma t^{-1}\left(1+O(\delta^2)\right),\\
\mu &= \eps^{-1/2} t^{3/2}\left(1+O(\delta^2)\right),
\end{align*}
и подставим сначала в \eqref{curve_delta}:
\begin{gather*}
  \delta\mu = \pm\ln c\mu^{1/2}+O(\mu^{-1/2}),\\
  \gamma\eps^{-1/2} t^{1/2}\left(1+O(\delta^2)\right)=
  \pm\ln c\eps^{-1/4} t^{3/4}+O(\delta^2)+O(\mu^{-1/2}),\\
  \gamma\eps^{-1/2} t^{1/2}=\pm\ln c\eps^{-1/4} t^{3/4}+O(\mu^{-1/2}),\\
  \gamma^\pm=\pm\eps^{1/2} t^{-1/2}\left(\ln c\eps^{-1/4} t^{3/4}+O(\eps^{1/4} t^{-3/4})\right),\\
\end{gather*}
а затем в \eqref{eigen_mu}:
\begin{gather*}
  \eps^{-1/2} t^{3/2}\left(1+O(\delta^2)\right)=3\pi(k-\frac{1}{8}\mp k_0(\lambda))
  +O(\mu^{-1/2}),\\
  \eps^{-1/2} t^{3/2}=3\pi(k-\frac{1}{8}\mp k_0(\lambda))
  +O(\mu^{-1/2}),\\
  \eps^{-1/2} \left(t_k^\pm\right)^{3/2}=3\pi(k-\frac{1}{8}\mp k_0(\lambda))
  +O(k^{-1/2}).
\end{gather*}
Остается заметить, что функция распределения собственных значений на кривых
$\gamma^\pm$ получается из формулы
\begin{equation*}
  (\eps^{-1/3}t_k^\pm)^{3/2}=3\pi k+O(1)
\end{equation*}
точно так же, как и при доказательстве теоремы \ref{th_orsommod_Nlam}.
\end{Proof}

\section{Поведение собственных значений в окрестности луча $[\node,-i\infty)$}
\label{par34}

Пусть параметр $\lambda$ лежит в области $D_\eps$, введенной при изучении
модельной задачи (параграф \ref{par31}). Тогда при $\sigma\to 0$
\begin{equation}
  |\arg\xi_j^\pm|\LE\pi-\frac{\ln{|\xi_j^\pm|^\varkappa}}{|\xi_j^\pm|^{3/2}},\quad
  \varkappa>\frac{3}{4}, \label{orsom_cond3a}
\end{equation}
что обеспечивает применимость леммы \ref{lm_aias}. Положим $j_1=0$, $j_2=2$ в
характеристическом определителе \eqref{detuxi}, а число $\xi_j^0$ в формуле
\eqref{ujxij} будем выбирать так, чтобы
$$
  \frac{2\pi}{3}j+\frac{\pi}{3}\LE
  \arg\xi_j^-\LE
  \arg\xi_j^0=
  \frac{2\pi}{3}(j+1)\LE
  \arg\xi_j^+\LE
  \frac{2\pi}{3}(j+1)+\frac{\pi}{3}.
$$

\begin{Lemma}
\label{lm_uxi_ray}
 При некотором выборе чисел $\xi_j^0$ интегралы \eqref{ujxij} имеют
 следующую асимптотику:
\begin{equation}
  u(\xi)=-\frac{e^{-\frac{2}{3}\xi^{3/2}}}{2\sqrt{\pi}\xi^{3/4}}
  (g(\sigma\xi)+O(\lambda\xi^{-3/2})),\quad
  \sigma\to 0.
  \label{prop3}
\end{equation}
Здесь $\xi=\xi_j^\pm$, а в качестве $u$, $g$ можно взять $u_j$, $g_j$ или
$u_j^*$, $g_j^*$.
\end{Lemma}
\begin{Proof}
Если $\lambda$ изменяется в ограниченной области, то при $j=0,1$ можно взять
$\xi_j^0=0$ и применить лемму \ref{lm_uxi_prop}. В противном случае эту лемму
непосредственно применить нельзя, но можно повторить с незначительными
изменениями вычисление интеграла по дуге окружности $|z|=|\xi|$. Для этого
надо подобрать $\xi_j^0$ так, чтобы
$$
  \frac{\xi_j}{\xi_j^0}=O(1),\qquad
  \RE\eta>0,\quad
  \frac{\eta}{\RE\eta}=O(1),\quad
  \eta=\frac{2}{3}((\xi_j^0)^{3/2}-\xi_j^{3/2}).
$$
При $j=0,1$ и достаточно больших $\lambda$ имеем:
$$
  \arg\xi_j^0=\frac{2\pi}{3}(j+1),\qquad
  |\arg\xi_j-\arg\xi_j^0|<\frac{\pi}{6}.
$$
Поэтому $\RE\xi_j^{3/2}\LE -\frac{|\xi_j|^{3/2}}{\sqrt{2}}$ и, полагая
$|\xi_j^0|=\frac{1}{2}\min\{|\xi_j^-|,|\xi_j^+|\}\LE\frac{1}{2}|\xi_j|$,
получим:
\begin{gather*}
  \RE\eta \GE
  \frac{2}{3}\left(-|\xi_j^0|^{3/2}+\frac{|\xi_j|^{3/2}}{\sqrt{2}}\right) \GE
  \frac{2}{3}\left(-\frac{|\xi_j|^{3/2}}{2\sqrt{2}}+
    \frac{|\xi_j|^{3/2}}{\sqrt{2}}\right)=
  \frac{1}{3\sqrt{2}}|\xi_j|^{3/2},\\
  |\eta|\LE \frac{2}{3}(|\xi_j^0|^{3/2}+|\xi_j|^{3/2}) \LE
  \frac{2}{3}|\xi_j|^{3/2}\left(\frac{1}{2\sqrt{2}}+1\right) \LE
  |\xi_j|^{3/2} \LE
  3\sqrt{2}\RE{\eta}.
   \\
  \frac{|\xi_j|}{|\xi_j^0|}=
  2\frac{|\xi_j|}{\min\{|\xi_j^-|,|\xi_j^+|\}}=
  2\frac{|\pm 1-\lambda|}{\min\{|-1-\lambda|,|1-\lambda|\}}=
  O(1).
\end{gather*}
При $j=2$ имеем:
$$
  \arg\xi_j^0=0,\qquad
  |\arg\xi_j|<\frac{\pi}{3}.
$$
Поэтому $\RE\xi_j^{3/2}\LE|\xi_j|^{3/2}$ и, полагая
$|\xi_j^0|=2\max\{|\xi_j^-|,|\xi_j^+|\}\GE 2|\xi_j|$, получим:
\begin{gather*}
  \RE\eta\GE
  \frac{2}{3}\left|\,|\xi_j^0|^{3/2}-|\xi_j|^{3/2}\right| \GE
  \frac{2}{3}|\xi_j^0|^{3/2}\left(1-\frac{1}{2\sqrt{2}}\right) \GE
  \frac{1}{3}|\xi_j^0|^{3/2},\\
  |\eta|\LE
  \frac{2}{3}\left(|\xi_j^0|^{3/2}+|\xi_j|^{3/2}\right) \LE
  \frac{2}{3}|\xi_j^0|^{3/2}\left(1+\frac{1}{2\sqrt{2}}\right) \LE
  |\xi_j^0|^{3/2} \LE
  3\RE\eta, \\
  \frac{|\xi_j|}{|\xi_j^0|} \LE
  \frac{1}{2}=
  O(1).
\end{gather*}
\end{Proof}

\begin{Lemma}
\label{lm_char_ray}
  Если $\lambda\in D_\eps$, то характеристическое уравнение
  приводится к виду
\begin{equation}
  \label{char_down}
  e^{\frac{4}{3}((\xi_0^+)^{3/2}-(\xi_0^-)^{3/2})}=
  1+|\lambda|^{-1/2}O(\sigma^{3/2}).
\end{equation}
\end{Lemma}
\begin{Proof}
Характеристическое уравнение имеет вид:
$$
\begin{vmatrix}
  u_0(\xi_0^+)-u_0(\xi_0^-) &
  u_2(\xi_2^+)-u_2(\xi_2^-) \\
  u_0^*(\xi_0^+)-u_0^*(\xi_0^-) &
  u_2^*(\xi_2^+)-u_2^*(\xi_2^-)
\end{vmatrix}
=0.
$$
Применим лемму \ref{lm_uxi_ray}, учитывая, что
\begin{gather*}
  g_j(\sigma\xi_j^\pm)=\sh\left[\alpha(1\mp 1)\right],\quad
  g_j^*(\sigma\xi_j^\pm)=\ch\left[\alpha(1\mp 1)\right],\\
  (\xi_2^\pm)^{3/2}=-(\xi_0^\pm)^{3/2},\quad
  (\xi_2^\pm)^{3/4}=-i(\xi_0^\pm)^{3/4}.
\end{gather*}
Получим:
$$
\begin{vmatrix}
  P_-E_-^{-1}[\sh 2\alpha]-P_+E_+^{-1}[0] &
  P_-E_-[i\sh 2\alpha]-P_+E_+[0] \\
  P_-E_-^{-1}[\ch 2\alpha]-P_+E_+^{-1}[1] &
  P_-E_-[i\ch 2\alpha]-P_+E_+[i]
\end{vmatrix}
=0,
$$
где
$$
  P_\pm=\frac{1}{2\sqrt{\pi}(\xi_0^\pm)^{3/4}},\qquad
  E_\pm=e^{\frac{2}{3}(\xi_0^\pm)^{3/2}},\qquad
  [c]=c+|\lambda|^{-1/2}O(\sigma^{3/2}).
$$
Преобразуем:
$$
\begin{vmatrix}
  E[\sh 2\alpha]-P[0] &
  [i\sh 2\alpha]-EP[0] \\
  E[\ch 2\alpha]-P[1] &
  [i\ch 2\alpha]-EP[i]
\end{vmatrix}
=0,
$$
где
$$
  P=P_+P_-^{-1}=\frac{(\xi_0^-)^{3/4}}{(\xi_0^+)^{3/4}},\qquad
  E=E_+E_-^{-1}=e^{\frac{2}{3}((\xi_0^+)^{3/2}-(\xi_0^-)^{3/2})}.
$$
Вычисляя определитель, получим:
\begin{gather*}
  -PE^2[i\sh 2\alpha]+E[i\sh 2\alpha \ch 2\alpha-i\sh 2\alpha \ch 2\alpha]+P[i\sh 2\alpha]=0, \\
  E^2[1]-E[0]-[1]=0, \\
  E^2=[1].
\end{gather*}
\end{Proof}

Сравним характеристическое уравнение \eqref{char_down} с соответствующим
уравнением \eqref{eq_char_model} для модельной задачи. Поскольку используемые
в параграфе \ref{par31} переменные $\xi_1$, $\xi_2$ равны соответственно
$\xi_0^-$, $\xi_0^+$, то \eqref{eq_char_model} можно переписать так:
\begin{equation*}
  e^{\frac{4}{3}((\xi_0^+)^{3/2}-(\xi_0^-)^{3/2})}=
  1+|\lambda|^{-3/2}O(\sigma^{3/2}).
\end{equation*}
Это уравнение отличается от \eqref{char_down} лишь бесконечно малым
слагаемым, поэтому справедлив следующий аналог теоремы \ref{th_orsommod_ray}
для задачи Орра--Зоммерфельда.
\begin{Theorem}
\label{th_orsom_ray} При достаточно малых $\eps>0$ все собственные значения
задачи \eqref{orsom1}, \eqref{orsom2}, лежащие в области $D_\eps$, являются
простыми и образуют монотонную последовательность чисел на отрицательной
мнимой оси. При этом на интервале
$\left(d_\eps,-i(\rho_{k_0}+O(\eps))\right)$ может находиться не более двух
собственных значений, а все остальные собственные значения $\lambda_k$ задачи
\eqref{orsom1}, \eqref{orsom2} в области $D_\eps$ имеют вид
$$
\lambda_k=-i(\rho_k+O(\eps)),\quad k=k_0+1,\;k_0+2,\;\ldots.
$$
Функция распределения собственных значений на луче $(\node,-i\infty)$
совпадает с соответствующей функцией распределения для модельной задачи
$$
N(\lambda,\eps)=\frac{\eps^{-1/2}}{\pi}
\int\limits_{-1}^{1}e^{-i\frac{\pi}{4}}\sqrt{x-\lambda}\,dx +\OB{1}.
$$
\end{Theorem}

\end{document}